\documentclass[final]{siamltex}
\usepackage[dvips]{graphicx}
\usepackage{amsmath}
\usepackage{latexsym}
\usepackage{amssymb}
\usepackage{algorithm,algorithmic}
\usepackage{color}
 \usepackage{mathptmx}

\title{Computational Krylov-based methods for   large-scale differential Sylvester matrix problems}
\author{ M. Hached \thanks{Laboratoire P.  Painlev\'e UMR 8524, UFR de Math\'{e}matiques, Universit\'e des Sciences et Technologies de Lille, IUT A, Rue de la Recherche, BP 179, 59653 Villeneuve d'Ascq Cedex, France; {\tt email:  mustapha.hached@univ-lille1.fr}} \and K. Jbilou\thanks{L.M.P.A, Universit\'e du Littoral C\^ote d'Opale,
50 rue F. Buisson BP. 699, F-62228 Calais Cedex, France. {\tt E-mail: ;
jbilou@univ-littoral.fr}.} }
\date{ }

\newcommand{\K}{\mathbb{K}}

\begin {document}
\maketitle
\begin{abstract}
In the present paper, we propose  Krylov-based methods for solving  large-scale  differential Sylvester matrix  equations having a low rank constant term. We present two new approaches for  solving such differential matrix equations. The first approach is based on the integral expression of the exact solution and  a Krylov method for  the computation of the exponential of a matrix times a block of vectors.   In the second approach, we first  project the initial problem onto a block (or extended block)  Krylov subspace and get a low-dimensional differential Sylvester  matrix  equation. The latter problem is then solved by some integration numerical methods such as BDF or Rosenbrock method and the obtained solution is used to build the low rank approximate solution of the original problem. We give some new  theoretical results such as a simple expression of the residual norm and upper bounds for the norm of the error. Some numerical experiments are given in order  to compare the two approaches.
\end{abstract}
\begin{keywords}
Extended block Krylov; Low rank; Differential Sylvester equations.
\end{keywords}
\begin{AMS}
65F10, 65F30
\end{AMS}
\pagestyle{myheadings}
\thispagestyle{plain}
\markboth{M. Hached and   K. Jbilou}{Computational Krylov-based methods ...}
%
%\begin{AMS}
%{\bf AMS Classification}: 65F10, 65F30
%\end{AMS}
%

%\date{ }
%
%\newtheorem{Remark}{Remark}
%\newcommand{\R}{\mathbb{R}}
%\newcommand{\K}{\mathbb{K}}
%\newcommand{\T}{\mathbb{T}}
%\newcommand{\M}{\mathbb{M}}
%\newcommand{\s}{\mathbb{S}}

%\textwidth 16.5cm
%%\textheight 20cm
%\oddsidemargin 0cm
%

%
%\maketitle
%

%
\section{Introduction}
In the present paper, we consider the  differential Sylvester 
 matrix  equation (DSE in short) of the form
\begin{equation}\label{sylv1}
\left\{
\begin{array}{l}
\dot X(t)=A(t)\,X(t)+X(t)\,B(t)+E(t)F(t)^T;\; (DSE) \\
 \;X(t_0)=X_0,\; \; t \in [t_0, \, T_f],
\end{array}
\right.
\end{equation}

\noindent where $ A(t) \in \mathbb{R}^{n \times n} $,   $ B(t) \in \mathbb{R}^{p \times p}$ and 
$E(t) \in \mathbb{R}^{n \times s} $ and  $F(t) \in \mathbb{R}^{p \times s} $ are   full rank matrices, with $ s\ll n,\,p $.  The initial condition is given in a factored form as $X_0=Z_0 \widetilde{Z} _0^T$ and the matrices $A$ and $B$ are assumed to be large and sparse.\\
Differential Sylvester equations play a fundamental role in many areas such as control,
filter design theory, model reduction problems, differential
equations and robust control problems \cite{abou03,corless}.
%For those  applications, the matrix  $A$ is generally sparse 
%and very large. 
For such differential matrix equations, only  a few attempts have been
made  for large problems. \\

Let us first recall the following theoretical result which gives an expression of the exact solution of \eqref{sylv1}.
\begin{theorem}\label{theo1}
\cite{abou03}
The unique solution of the  general differential Sylvester  equation
\begin{equation}
\displaystyle {\dot X}(t)=A(t)\,X(t)+X(t)\,B(t)+M(t);\;\; X(t_0)=X_0
\end{equation}
is defined by
\begin{equation}\label{solexacte1}
X(t) = \Phi_A(t,0)X_0{\Phi_{B^T}^T(t,t_0)}+\int_{t_0}^t \Phi_A(t,\tau)M(\tau){\Phi_{B^T}^T(t,\tau)}d\tau.
\end{equation}
where the transition matrix $\Phi_A(t,t_0)$ is the unique solution to the problem
$$\displaystyle {\dot \Phi}_A(t,t_0)=A(t) \Phi_A(t,t_0),\;\; \Phi_A(t_0,t_0)=I.$$
Futhermore, if $A$ is assumed to be a constant matrix, then we have
\begin{equation}\label{solexacte2}
X(t)=e^{(t-t_0)A}X_0e^{(t-t_0)B}+\int_{t_0}^t e^{(t-\tau)A}M(\tau)e^{(t-\tau)B}d\tau.
\end{equation}
\end{theorem}
We notice that the problem \eqref{sylv1} is equivalent to the linear ordinary differential equation
\begin{equation}\label{kron}
\left\{
\begin{array}{c c l}
\dot{x}(t)& =& \mathcal{A}(t)x(t)+b(t) \\
x_0 & = & vec(X_0)
\end{array}
\right.
\end{equation}
where $\mathcal{A}= I \otimes A(t) + B^T(t) \otimes I$, $x(t)=vec(X(t)$ and $b(t) = vec(E(t)F(t)^T)$, where $vec(Z)$ is the long  vector obtained by stacking the columns of the matrix $Z$, forming a sole column. For moderate size problems, it is then possible to directly apply an integration method to solve  \eqref{kron}. However, this approach is not suitable  for large problems. From now on, we assume that the matrices $A$ and $B$ are time independent.\\

 In the present paper, we will consider projection methods onto extended block Krylov   (or block Krylov) subspaces  associated to the pairs $(A,E)$ and $(B^T,F)$ defined as follows
$$\K_m(A,E)={\rm range}(E,AE,\ldots,A^{m-1}E)$$
for block Krylov subspaces, or
$${\cal K}_m(A,E)={\rm range}(A^{-m}E,\ldots,A^{-1}E,E,AE,\ldots,A^{m-1}E)$$
for extended block Krylov subspaces when the matrix $A$ is nonsingular. Notice that  the extended block Krylov subspace ${\cal K}_m(A,E)$ is a sum of two block Krylov subspaces associated to the pairs $(A,E)$ and $(A^{-1},A^{-1}E)$:
$$ {\cal K}_m(A,E)=\K_{m}(A,E)\, + \, \K_m(A^{-1},A^{-1}E). $$
To compute an orthonormal basis $\{V_1,\ldots,V_m\}$,  where $V_i$ is of dimension $n\times d$ where $d=s$ for the block Krylov and $d=2s$ in  the extended  block Krylov case, two  algorithms have been defined: the first one is the well known block Arnoldi algorithm and the second one is the extended block Arnoldi algorithm \cite{druskin98,simoncini1}; see Appendix A for the description of both algorithms. \\These  algorithms generate the blocks  $ V_1,V_2,\ldots,V_m  $, $ V_ i \in \mathbb{R}^{n \times d} $ such that  their columns form an orthonormal basis of the block Krylov subspace $\K_{m}(A,E)$ (with $d=s$) or   the extended block Arnoldi ${\cal K}_m(A,E)$ (with $d=2s$). \\ Both algorithms compute also $d \times d$ block upper Hessenberg matrices $ {\mathcal T}_{m,A} = {\mathcal V}_m^T\,A\,{\mathcal V}_m $. The following algebraic relations are satisfied
\begin{eqnarray}
A\,{\cal V}_m & = & {\cal V}_{m+1}\,{\widehat {\cal T}}_{m,A}, \\
              & = & {\cal V}_m\,{\cal T}_{m,A} + V_{m+1}\,T_{m+1,m}\,\widetilde{E}_m^T,
\end{eqnarray}
where  ${{\widehat {\cal T}}_{m,A}} = {\cal V}_{m+1}^T\,A\,{\cal V}_m ; $  $ T_{i,j} $ is the  $ (i,j) $ block of $  {\widehat {\cal T}_{m,A}} $ of size $ d\times d $, and $ \widetilde{E}_m = [ O_{d \times (m-1)d}, I_{d} ]^T $ is the matrix formed with the last $ d$ columns of the $ md \times md$ identity matrix $ I_{md} $ where $d=s$ for the block Arnoldi and $d=2s$ for the extended block Arnoldi.\\ When the matrix $A$ is nonsingular and when the computation of the products $W=A^{-1}V$ is not difficult (which is the case for sparse and structured matrices), the use of the extended block Arnoldi is to be preferred. \\

The paper is organized as follows: In Section 2, we present a first approach based on the approximation of the exponential of a matrix times a block using a Krylov projection method. We give some theoretical results such as a simple expression of the norm of the residual and upper bounds for the norm of the error and perturbation results.  In Section 3,   the initial differential Sylvester matrix equation is projected onto a block (or extended block) Krylov subspace. The obtained low dimensional differential Sylvester  equation is solved by using the well known Backward Differentiation Formula (BDF) and Rosenbrock methods. The last section is devoted to some numerical experiments.
\\
Throught the paper,  $\Vert . \Vert$  and $\Vert \,.\, \Vert_F$ will denote the 2-norm and the Frobenius norm, respectively.  

%******************************************************************************

\section{Solutions via  of the matrix exponential approximation }

In this section, we give a new  approach for computing approximate solutions to large differential Sylvester equations \eqref{sylv1}. \\ 
We recall that the exact solution to \eqref{sylv1} can be expressed as follows 
\begin{equation}\label{solexacte3}
X(t)=e^{(t-t_0)A}X_0e^{(t-t_0)B}+\int_{t_0}^t  e^{(t-\tau)A}\, EF^Te^{(t-\tau)B}\, d\tau.
\end{equation}
For our first  approach, we use this expression of $X(t)$ to obtain low rank approximate solutions. We first  approximate  the factors $ e^{(t-\tau)A}E$ and   $ e^{(t-\tau)B^T}F$ and   then, use  a quadrature method to   compute the desired approximate solution. As the matrices $e^{(t-\tau)A}$ and  $ e^{(t-\tau)B^T}$  are  large and  could be dense even though $A$ and $B$  are sparse, computing the exponential is not recommended. However, in our problem, the computation of $e^{(t-\tau)A}$ and $ e^{(t-\tau)B^T}$  are not needed explicitly as we will rather consider the products $e^{(t-\tau)A}\, E$ and $ e^{(t-\tau)B^T}F$  for which approximations via projection methods onto block or extended block Krylov subspaces are well suited. 
\\
 In what follows, we consider projections onto extended block Krylov (or just block Krylov) subspaces.  
Let $\mathcal{V}_m=[V_1,\ldots,V_m]$ and  $\mathcal{W}_m=[W_1,\ldots,W_m]$ be the orthogonal  matrices whose  columns form an orthonormal  basis of the subspace ${\cal K}_m(A,E)$ and ${\cal K}_m(B^T,F)$, respectively.    
Following \cite{saad1,saad2,vorst1}, an approximation to $Z_A=e^{(t-\tau)A}\, E$ can be obtained as
\begin{equation}
\label{expA}
Z_{m,A}(\tau) = \mathcal{V}_m e^{(t-\tau)\mathcal{T}_{m,A}}\, \mathcal{V}_m^T E
\end{equation}
where $\mathcal{T}_{m,A}=\mathcal{V}^T_{m} A \mathcal{V}_m$. In the same way, an approximation to $ e^{(t-\tau)B^T}F$, is given by
\begin{equation}
\label{expB}
Z_{m,B}(\tau) = \mathcal{W}_m e^{(t-\tau)\mathcal{T}_{m,B}}\, \mathcal{W}_m^T F,
\end{equation}
where $\mathcal{T}_{m,B}=\mathcal{W}^T_{m} B^T \mathcal{W}_m$. 
Therefore, the integrand in the expression \eqref{solexacte3} can be approximated as 
\begin{equation}
\label{exp2}
e^{(t-\tau)A}EF^Te^{(t-\tau)B} \approx Z_{m,A}(\tau) Z_{m,B}(\tau)^T.
\end{equation}
If for simplicity, we assume that $X(0)=0$, an approximation to the solution of the differential Sylvester  equation \eqref{sylv1} can be expressed as the product
\begin{equation}
\label{exp3}
X_m(t) = \mathcal{V}_m G_m(t) {\mathcal{W}_m}^T,\; t \in [t_0,\,T_f],
\end{equation}
where 
\begin{equation}
\label{gm}
G_m(t) = \displaystyle  \int_{t_0}^t  Z_{m,A}(\tau) Z_{m,B}^T(\tau) d\tau,
\end{equation}
with   $E_m=\mathcal{V}^T_{m} E$ and $F_m=\mathcal{W}^T_{m} F$. \\

\noindent The next result shows that the $md \times md$ matrix function $G_{m}$ is  solution of a low-order differential Sylvester  matrix equation.
\begin{proposition}
Let $G_m(t)$ be the matrix function  defined by \eqref{gm}, then it satisfies the following low-order differential Sylvester  matrix equation
\begin{equation}\label{low2}
{\dot G}_m(t) = \mathcal{T}_{m,A} G_m(t) + G_m(t){\mathcal{T}_{m,B}}^T  +E_mF_m^T,\;\;  t \in [t_0,\,T_f].
\end{equation}
\end{proposition}
\medskip
{\bf Proof.}  
The proof can be easily derived from the expression \eqref{gm} and the result of Theorem \ref{theo1}.

\medskip
\noindent As a consequence, introducing  the residual $ R_m(t) = \displaystyle {\dot X}_m(t)-A\,X_m(t)-X_m(t)\,B- EF^T $ associated to the approximation $X_m(t)$,  we have the following relation
\begin{eqnarray*}
{\cal V}^T_{m} R_m(t) {\cal W}_m &= & {\cal V}^T_m ({\dot X}_m(t) -AX_m(t)-X_m(t)B-EF^T) {\cal W}_m\\
& = &  {\dot G}_m(t) -\mathcal{T}_{m,A} G_m(t) - G_m(t){\mathcal{T}_{m,B}}^T  -E_mF_m^T\\
& = & 0,
\end{eqnarray*}
which shows that the residual satisfies a Petrov-Galerkin condition.\\

\noindent As mentioned earlier and for our first exponential-based approach,  once  $ Z_{m,A}(\tau)$ and  $ Z_{m,B}(\tau)$ are  computed, we use a quadrature method to approximate the integral \eqref{gm} in order to get an approximation of  $G_m(t)$ and hence to compute $X_m(t)$ from \eqref{exp3}.

\noindent The computation of $ X_m(t) $
(and of $ R_m(t) $) becomes expensive as $ m $ increases. So, in
order to stop the iterations, one has to test if $ \parallel R_m(t)
\parallel < \epsilon $ without having to compute extra products
involving the matrices $ A $ and $B$. The next result shows how to compute
the  norm of $ R_m(t) $ without forming the approximation $
X_m(t) $ which is computed in a factored form only when convergence
is achieved. 
\medskip
\begin{proposition} \label{t2}
Let $ X_m(t) = {\cal V}_mG_m(t){\cal W}_m^T $ be the approximation obtained at step $ m $ by the  block (or extended block) Arnoldi  method. Then the residual $ R_m(t) $ satisfies the relation

\begin{equation}
\label{result2}
\parallel R_m(t) \parallel_F^2 =  \parallel T_{m+1,m}^A {\bar G}_{m}(t) \parallel_F^2 + \parallel T_{m+1,m}^B {\bar G}_{m}(t) \parallel_F^2,
\end{equation}
and for the 2-norm, we have
\begin{equation}
\label{resultnorm2}
\parallel R_m(t) \parallel = \displaystyle \max \{ \parallel T_{m+1,m}^A {\bar G}_{m}(t) \parallel,  \parallel T_{m+1,m}^B {\bar G}_{m}(t) \parallel\},
\end{equation}
where $ {\bar G}_m $ is the $ d \times md $  matrix corresponding to the last $ d$ rows of $ G_m $ where $d=2s$ when using the extended block Arnoldi algorithm and $d=s$ with the block Arnoldi algorithm. 
\end{proposition} 
\medskip
{\bf Proof.}  . 
The proof  comes from the fact that the residual $R_m(t)$ can be expressed as 
\begin{equation}
\label{rm}
R_m(t) = {\cal V}_{m+1} \left ( 
\begin{array}{cc}
 {\dot G}_m(t) -\mathcal{T}_{m,A} G_m(t) - G_m(t){\mathcal{T}_{m,B}}^T  -E_mF_m^T & -T_{m+1,m}^B {\bar G}_{m}(t)\\
  T_{m+1,m}^A {\bar G}_{m}(t)  & 0 
\end{array}
\right )\, {\cal W}_{m+1}^T,
\end{equation}
where  $G_m(t)$ solves the low dimensional problem \eqref{low2}. Therefore, we get
\begin{eqnarray*}
\Vert R_m(t) \Vert_F^2 & = & \left  \Vert \left ( 
\begin{array}{cc}
0 & -T_{m+1,m}^B {\bar G}_{m}(t)\\
  T_{m+1,m}^A {\bar G}_{m}(t) & 0 
\end{array}
\right ) \right \Vert_F^2\\
& = & \parallel T_{m+1,m}^A {\bar G}_{m}(t) \parallel_F^2 + \parallel T_{m+1,m}^B {\bar G}_{m}(t) \parallel_F^2.
\end{eqnarray*}

\noindent To prove the expression  \eqref{resultnorm2} with the 2-norm , let us first remark that if  $$M=\left (  \begin{array}{cc} 
0 & M_1\\
M_2 & 0
\end{array}\right),\;\;  {\rm then}\; \; M^TM=\left (  \begin{array}{cc} 
M_1^T M_1 & 0\\
0 & M_2 ^T M_2 
\end{array}\right),$$
which shows that the singular values  of $M$ are the sum of the singular values of $M_1$ and those of $M_2$ which implies that$$ \Vert M \Vert= \sigma_{max}(M)= \displaystyle \max\{ \sigma_{max}(M_1), \sigma_{max}(M_2)\}= \max\{\Vert M_1 \Vert, \Vert M_2 \Vert\}.$$
Therefore, using this remark and the fact that  
$$\Vert R_m(t) \Vert = \left  \Vert \left ( 
\begin{array}{cc}
0 & -T_{m+1,m}^B {\bar G}_{m}(t)\\
  T_{m+1,m}^A {\bar G}_{m}(t) & 0 
\end{array}
\right ) \right \Vert,$$
the result follows.\\

\medskip
\noindent 
%The result given by \eqref{result2} or by \eqref{resultnorm2} is very  important in practice, as it allows us to stop the iterations when convergence is achieved without computing the approximate solution $X_m(t)$. \\   

The approximate  solution
$X_m(t)$   is computed only when convergence is achieved and in a factored form  which is very important for storage requirements in large-scale problems. This procedure is described as follows.\\
Consider the singular value decomposition of the  
matrix $G_m(t)=U\, \Sigma\, V$ where $\Sigma$ is the
diagonal matrix of the  singular values of $G_m(t)$ sorted in
decreasing order. Let $U_l$ be the $md \times l$ matrix  of   the first $l$ columns of  $U$ 
corresponding to the $l$ singular values  of magnitude greater than
 some tolerance $dtol$. We obtain the
truncated singular value decomposition  $G_m(t) \approx U_l\, \Sigma_l\,
V_l^T$ where $\Sigma_l = {\rm diag}[\lambda_1, \ldots, \lambda_l]$.
 Setting ${\widetilde Z}_{m,A}(t)={\cal V}_m \, U_l\, \Sigma_l^{1/2}$ and ${\widetilde Z}_{m,B}(t)={\cal W}_m \, V_l\, \Sigma_l^{1/2}$ , it
follows that
\begin{equation}
\label{approx}
X_m(t) \approx {\widetilde Z}_{m,A}(t) {\widetilde Z}_{m,B}(t)^T.
\end{equation}
Therefore, only the matrices ${\widetilde Z}_{m,A}(t)$ and ${\widetilde Z}_{m,B}(t)$ are needed.\\ 

\noindent The following result  shows  that the approximation $X_m $ is an exact solution of a perturbed  differential  Sylvester   equation.\\

\begin{proposition}\label{ppertu}
Let $X_m(t)$ be the approximate solution given by \eqref{exp3}. Then we have 
\begin{equation}
\label{pertu}
\displaystyle {\dot X}_m(t)=(A-F_{m,A})\,X_m(t)+X_m(t)\,(B-{F}_{m,B})+EF^T.
\end{equation}
where $ F_{m,A} = V_{m+1}\,T_{m+1,m}^A\,V_{m}^T $ and $F_{m,B}=W_m(T_{m+1,m}^B)^T W_{m+1}^T$
\end{proposition}

\medskip

{\bf Proof.}  . 
As $X_m(t)={\cal V}_m G_m(t) {\cal W}_m^T$, we have 
\begin{equation}
\label{eq1}
{\dot X}_m(t)-(AX_m(t)+X_m(t)B+EF^T)=  {\cal V}_m {\dot G}_m(t) {\cal W}_m^T-(A{\cal V}_m G_m(t) {\cal W}_m^T +B{\cal V}_m G_m(t) {\cal W}_m^T +EF^T).
\end{equation}
Now, using the fact that $$A\,{\cal V}_m  =  {\cal V}_{m}\, {\cal T}_{m,A} + V_{m+1}T_{m+1,m}^A {\widetilde E}_m^T, \; {\rm and} \; B^T\,{\cal W}_m  =  {\cal W}_{m}\,{\cal T}_{m,B}+W_{m+1}T_{m+1,m}^B {\widetilde E}_m^T,$$ equation \eqref{eq1} becomes
\begin{eqnarray*}
\label{eq2}
{\dot X}_m(t)-(AX_m(t)+X_m(t) B+EF^T) & = & {\cal V}_m \dot G_m(t) {\cal W}_m^T-( [{\cal V}_{m}\, {\cal T}_{m,A} + V_{m+1}T_{m+1,m} {\widetilde E}_m^T]G_m(t) {\cal W}_m^T\\ & + & {\cal V}_mG_m(t)  [ {\cal W}_{m}\,{\cal T}_{m,B}+ W_{m+1}T_{m+1,m}^B {{\widetilde E}_m}^T]^T +EF^T ).
\end{eqnarray*}
Therefore
\begin{eqnarray*}
\label{eq3}
{\dot X}_m(t)-(AX_m(t)+X_m(t) B+EF^T) & = & {\cal V}_m [ {\dot G}_m(t) - {\cal T}_{m,A}  G_m(t) -G_m(t) {\cal T}_{m,B}^T -EF^T] {\cal W}_m\\
& - & ( V_{m+1}T_{m+1,m}^A {\widetilde E}_m^T G_m(t) {\cal W}_m + {\cal V}_mG_m(t)  {\widetilde E}_m(T_{m+1,m}^B)^T W_{m+1}^T ).
\end{eqnarray*}
On the other hand we have ${\cal V}_mG_m(t)=X_m(t) {\cal W}_m$,  $G_m(t) {\cal W}_m^T= {\cal V}_m^TX_m(t)$, ${\cal V}_m{\widetilde E}_m=V_m$, ${\cal W}_m{\widetilde E}_m=W_m$ and $EF^T={\cal V}_m E_mF_m^T {\cal W}_m^T $. So using these relations  and the fact that $G_m$ solves the low dimensional differential Sylvester equation \eqref{low2}, we obtain the desired result.

\medskip

\noindent The next result states that the error  ${\cal E}_m(t)=X(t)-X_m(t)$ satisfies also a differential Sylvester  matrix equation.
\medskip
\begin{proposition}\label{err1}
Let $X(t)$ be the exact solution of \eqref{sylv1} and let $X_m(t)$ be the approximate solution obtained at step $m$.  The error ${\cal E}_m(t)=X(t)-X_m(t)$  satisfies the following equation
\begin{equation}\label{pertu2}
\displaystyle \dot {\cal E}_m(t)-A{\cal E}_m(t)-{\cal E}_m(t)B=F_{m,A}X_m(t)+X_m(t) F_{m,B}=-R_m(t),
\end{equation}
where $F_{m,A}$ and $F_{m,B}$ are defined in Proposition \ref{ppertu} and $R_m(t)=  \displaystyle \dot{X}_m(t)-A\,X_m(t)-X_m(t)\,B- EF^T$.
\end{proposition}

\medskip

{\bf Proof.}  
The result is easily obtained by subtracting the  equation \eqref{pertu} from the initial differential Sylvester  equation  \eqref{sylv1}.\\

\noindent Notice that from Proposition \ref{err1}, the error ${\cal E}_m(t)$ can be expressed in the integral form as follows
\begin{equation}
\label{error3}
 {\cal E}_m(t)=e^{(t-t_0)A}{ \cal E}_{m,0}e^{(t-t_0)B}-\int_{t_0}^t e^{(t-\tau)A}R_m(\tau)e^{(t-\tau)B}d\tau,\; t \in [t_0,\, T_f].
\end{equation}
where  $\mathcal{E}_{m,0}=\mathcal{E}_m(0)$.

\medskip
\noindent Next, we give  an upper bound for the norm of the error by using the 2-logarithmic norm defined by  $\mu_2(A)=\displaystyle \lim_{h \rightarrow 0^+} \frac{\Vert I + hA \Vert_2-1}{h}=\displaystyle \frac{1}{2} \lambda_{max}(A+A^T)$. \\

\begin{proposition}
\label{t3}
Assume that the matrices $A$ and $B$ are such that $\mu_2(A)+\mu_2(B) \ne 0$. Then at step $m$ of the extended block Arnoldi (or block Arnoldi) process, we have the following upper bound for the norm of the error  $\mathcal{E}_m(t) = X(t)-X_m(t)$,
\begin{equation}
\label{upperbound}
\parallel \mathcal{E}_m(t) \parallel  \le  \Vert \mathcal{E}_{m,0} \Vert  e^{(t-t_0)(\mu_2(A)+\mu_2(B))}+ \alpha_m \frac{e^{(t-t_0)(\mu_2(A)+\mu_2(B))}-1}{ \mu_2(A)+\mu_2(B)},\\
\end{equation}
where  $\alpha_m$ is given by  $\alpha_m=\displaystyle 
\max_{\tau \in [t_0,\, t]}   \left( \displaystyle \max \{ \parallel T_{m+1,m}^A {\bar G}_{m}(t) \parallel_2,  \parallel T_{m+1,m}^B {\bar G}_{m}(t) \parallel_2\} \right)$. The matrix  $ {\bar G}_m $ is the $ d \times md $  matrix corresponding to the last $ d $ rows of $ G_m $. 
\end{proposition}

\medskip

{\bf Proof.}  
We first point out that $\parallel e^{tA}  \parallel \le e^{\mu_2(A)t}$. Using the expression \eqref{error3} of $\mathcal{E}_m(t) $, we  obtain the following relation 
\begin{equation*}
\parallel \mathcal{E}_m(t) \parallel  \le   \Vert \displaystyle  e^{(t-t_0)A}\mathcal{E}_{m,0}e^{(t-t_0)B} \Vert+   \displaystyle \int_{t_0}^t \Vert e^{(t-\tau)A}R_m(\tau)e^{(t-\tau)B} \Vert d \tau.
\end{equation*}
Therefore, using \eqref{error3} and the fact that $\parallel e^{(t-\tau)A} \parallel \le e^{(t-\tau) \mu_2(A)}$, we get 
\begin{eqnarray*}
\parallel \mathcal{E}_m(t) \parallel  &  \le  &  \Vert {\cal E}_{m,0} \Vert  e^{(t-t_0)(\mu_2(A)+\mu_2(B))} + \displaystyle \max_{\tau \in [t_0,t]} \parallel  R_m(\tau ) \parallel \displaystyle \int_{t_0}^t e^{(t-\tau) \mu_2(A)} e^{(t-\tau) \mu_2(B)}  d\tau\\
 & = &\Vert \mathcal{E}_{m,0} \Vert  e^{(t-t_0)(\mu_2(A)+\mu_2(B))} +  \displaystyle \max_{\tau \in [t_0,t]} \parallel  R_m(\tau ) \parallel e^{t(\mu_2(A)+\mu_2(B))} \displaystyle \int_{t_0}^t e^{-\tau (\mu_2(A)+\mu_2(B))} d\tau. \\
% & \le &  \Vert \mathcal{E}_{m,0} \Vert  e^{(t-t_0)(\mu_2(A)+\mu_2(B))} + \displaystyle \max_{\tau  \in [t_0,t]} \parallel  R_m(\tau ) \parallel \, e^{t(\mu_2(A)+\mu_2(B))} \\
% & \times  & \frac{e^{-t_0(\mu_2(A)+\mu_2(B))}-e^{-t(\mu_2(A)+ \mu_2(B)}}{\mu_2(A)+\mu_2(B)}.
 \end{eqnarray*}
Hence
 \begin{equation}
 \label{equp1}
 \parallel \mathcal{E}_m(t) \parallel \le \Vert \mathcal{E}_{m,0} \Vert \displaystyle  e^{(t-t_0)(\mu_2(A)+\mu_2(B))} +  \displaystyle \max_{\tau  \in [t_0,t]} \parallel  R_m(\tau ) \parallel \displaystyle  \frac{e^{(t-t_0)(\mu_2(A)+\mu_2(B))}-1}{ \mu_2(A)+\mu_2(B)}.
 \end{equation}
Using the result of Proposition \ref{t2},  we obtain $\displaystyle \max_{\tau  \in [t_0,t]} \parallel  R_m(\tau ) \parallel=\alpha_m$ and then 
$$\parallel \mathcal{E}_m(t) \parallel \le \Vert \mathcal{E}_{m,0} \Vert  e^{(t-t_0)(\mu_2(A)+\mu_2(B))} +  \alpha_m \, \displaystyle  \frac{e^{(t-t_0)(\mu_2(A)+\mu_2(B))}-1}{ \mu_2(A)+\mu_2(B)}.$$

\medskip

\noindent Notice that if  the matrices $A$ and $B$ are stable (\textit{ie} all the eigenvalues are in the open half plane) then  $\mu_2(A) <0$ and $\mu_2(B)<0$ which ensures the condition of Proposition \ref{t3} is satisfied.  
Notice also that since   $R_m(\tau)=-F_{m,A}X_m(\tau)-X_m(\tau) F_{m,B}$, where $ F_{m,A} = V_{m+1}\,T_{m+1,m}^A\,V_{m}^T $ and $F_{m,B}=W_m(T_{m+1,m}^B)^T W_{m+1}^T$, we get 
$$\Vert R_m(\tau) \Vert  \le \displaystyle \max_{\tau \in [t_0,t]} \Vert {\bar G}_m(\tau)\Vert \, \left (\Vert T_{m+1,m}^A \Vert + \Vert T_{m+1,m}^B \Vert \right ).$$
Hence, replacing in \eqref{equp1}, we get the new upper bound 
\begin{equation}
 \label{equp2}
 \parallel \mathcal{E}_m(t) \parallel  \le \Vert \mathcal{E}_{m,0} \Vert  e^{(t-t_0)(\mu_2(A)+\mu_2(B))} + \beta_m \displaystyle  \frac{e^{(t-t_0)(\mu_2(A)+\mu_2(B))}-1}{ \mu_2(A)+\mu_2(B)},
 \end{equation}
where $$\beta_m= \displaystyle \max_{\tau \in [t_0,t]} \Vert {\bar G}_m(\tau)\Vert \, \left(\Vert T_{m+1,m}^A \Vert + \Vert T_{m+1,m}^B \Vert \right).$$

\begin{figure}[H]
	\begin{center}
		\includegraphics[width=12cm,height=7cm]{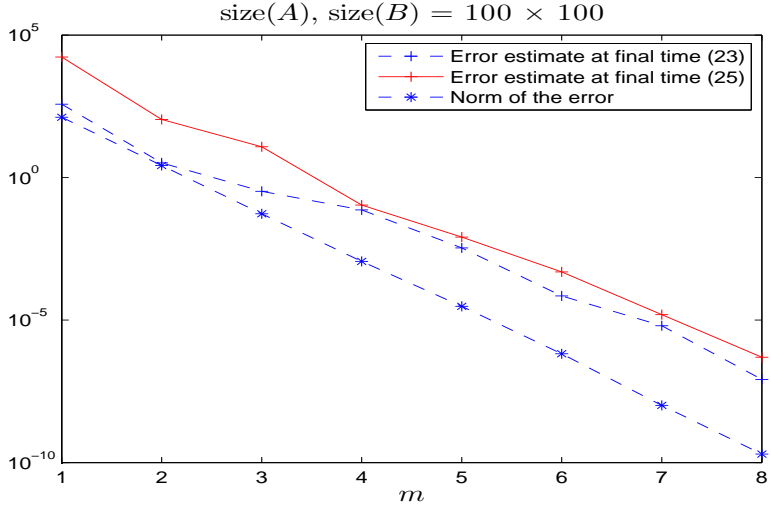}
		\caption{ Norm of the error \textit{vs} number of Arnoldi iterations $m$}\label{Figure1}
	\end{center}
\end{figure}

\noindent In Figure \ref{Figure1}, we compare the computed error to the two error upper bounds given by Formulae (\ref{upperbound}) and (\ref{equp2}) for $A$ and $B$ being two $100 \times 100$ matrices obtained by the finite differences discretization of linear differential operators on the unit square $[0,1]\times [0,1]$ with homogeneous Dirichlet boundary conditions. Matrices $E$ and $F$ were chosen as rank 2 matrices which entries are randomly generated over the interval $[0,1]$. In order to compute the error, we took the approximate solution given by the integral form of the solution as a reference.
\noindent We observe that the bound (\ref{upperbound}) stated in Proposition \ref{t3} is slightly better in this example.

\noindent Next, we give another  upper bound for the norm of the error ${\cal E}_m(t) $ .
\medskip
\begin{proposition}
\label{tr4}
Let $X(t)$ be the exact solution to \eqref{sylv1} and let $X_m(t)$ be the approximate solution obtained at step $m$.  Then we have
\begin{equation}
\label{err3}
\Vert {\cal E}_m(t) \Vert \le  \Vert F\,  \Vert e^{t\mu_2(B)} \, \Gamma_{1,m}(t) +  \Vert E_m\,  \Vert e^{t\mu_2(A)} \Gamma_{2,m}(t),
\end{equation}
where $$\Gamma_{1,m}(t) =\displaystyle \int_{t_0}^t e^{-\tau\mu_2(B)} \,  \Vert Z_A(\tau)-Z_{m,A}(\tau) \Vert  d \tau, \; \Gamma_{2,m}(t)=
\displaystyle \int_{t_0}^t e^{-\tau\mu_2(A)} \Vert  Z_B(\tau)-Z_{m,B}(\tau) \Vert   d \tau.$$
\end{proposition}
\medskip
{\bf Proof.}  
From the expressions of $X(t)$ and $X_m(t)$, we have 
\begin{equation}\label{errr3}
\Vert {\cal E}_m(t)\Vert  = \left  \Vert \int_{t_{0}}^t \left (Z_A(\tau)Z_B(\tau)^T -Z_{m,A}(\tau)Z_{m,B}(\tau)^T \right ) d \tau  \right \Vert,
\end{equation}
where $Z_{m,A}=  {\cal V}_m e^{(t-\tau)\mathcal{T}_{m,A}}E_m$, $ Z_{m,B}(\tau)= {\cal W}_m e^{(t-\tau)\mathcal{T}_{m,B}}F_m$, $Z_A(\tau)=e^{(t-\tau)A}E$ and $Z_B(\tau)=e^{(t-\tau)B}F$.  Then, using the relation
$$ Z_A(\tau)Z_B(\tau)^T -Z_{m,A}(\tau)Z_{m,B}(\tau)^T= (Z_A(\tau)-Z_{m,A}(\tau))Z_B^T+ Z_{m,A}(\tau) (Z_B(\tau)-Z_{m,B}(\tau))^T,$$ we obtain 
\begin{eqnarray*}
\Vert  Z_A(\tau)Z_B(\tau)^T -Z_{m,A}(\tau)Z_{m,B}(\tau)^T \Vert  & \le &   \Vert Z_B(\tau) \Vert\, \Vert (Z_A(\tau)-Z_{m,A}(\tau)) \Vert \\ &  + & \Vert Z_{m,A}(\tau) \Vert \, \Vert  (Z_B(\tau)-Z_{m,B}(\tau)) \Vert.
\end{eqnarray*}
Now as $\Vert Z_B(\tau) \Vert \le e^{(t-\tau)\mu_2(B)} \Vert F \Vert$ and since $\mu_2({\cal T}_{m,A}) \le \mu_2(A)$, we also have $\Vert Z_{m,A}(\tau) \Vert \le e^{(t-\tau)\mu_2({\cal T}_{m,A})} \Vert E_m \Vert \le e^{(t-\tau)\mu_2(A)} \Vert E_m \Vert$.
Using all these relations in \eqref{errr3}, we get
\begin{eqnarray*}
\label{err4}
\Vert {\cal E}_m(t) \Vert   & \le  & \int_{t_0}^t \left [e^{(t-\tau)\mu_2(B)} \Vert F\,  \Vert Z_A(\tau)-Z_{m,A}(\tau) \Vert + e^{(t-\tau)\mu_2(A)} \Vert E_m \Vert\,  \Vert  Z_B(\tau)-Z_{m,B}(\tau)\Vert   \right] d\tau\\
& \le & \Vert F\,  \Vert e^{t\mu_2(B)} \, \int_{t_0}^t e^{-\tau\mu_2(B)} \,  \Vert Z_A(\tau)-Z_{m,A}(\tau) \Vert  d\tau \\
&+ &  \Vert E_m\,  \Vert e^{t\mu_2(A)} \, \int_{t_0}^t   e^{-\tau\mu_2(A)}\Vert  (Z_B(\tau)-Z_{m,B}(\tau)) \Vert   d \tau,
\end{eqnarray*}
which ends the proof.\\

\noindent One can use some known results \cite{hochbruck,saad2} to derive  upper bounds for $\Vert Z_A(\tau)-Z_{m,A}(\tau) \Vert$ and
$\Vert  Z_B(\tau)-Z_{m,B}(\tau) \Vert $, when using Krylov or block Krylov subspaces.  For general matrices $A$ and $B$, we can use the following result to get upper bounds for $\Vert Z_A(\tau)-Z_{m,A}(\tau) \Vert$ and
$\Vert  Z_B(\tau)-Z_{m,B}(\tau) \Vert $.\\

\begin{proposition}
When using the extended block Arnoldi (or the block Arnoldi), we get the following upper bound for the exponential approximation error $ e_{m,A}(\tau)=Z_A(\tau)-Z_{m,A}(\tau) $:
\begin{equation}\label{erm}
\Vert e_{m,A}(\tau) \Vert \le \Vert T_{m+1,m}^A \Vert \int_0^{\tau}  e ^{(u-\tau)\nu(A)} \Vert L_{m,A}(u) \Vert d u,
\end{equation}
where $L_{m,A}(u)={\widetilde E}_me^{(t-u) {\cal T}_{m,A}}E_m$ and $\nu(A)=\lambda_{min} \left ( \displaystyle \frac{A+A^T}{2}\right)$.
\end{proposition}

\medskip

{\bf Proof.}  
We have $$Z_A(\tau) =  e^{(t-\tau) A} E,\;\;  {\rm and} \;\; Z_{m,A}(\tau)={\cal V}_m e^{(t-\tau) {\cal T}_{m,A}} E_m.$$
Then ${Z^{\prime}_A}(\tau) =  -Ae^{(t-\tau) A} E=-AZ_A(\tau)$, 
and
$${Z^{\prime}_{m,A}}(\tau)=-{\cal V}_m{\cal T}_{m,A} e^{(t-\tau) {\cal T}_{m,A}} E_m=-[A{\cal V}_m -V_{m+1} T_{m+1,m}^A {\widetilde E}_m]e^{(t-\tau) {\cal T}_{m,A}}E_m.$$
Hence, 
\begin{equation}
{Z^{\prime}_{m,A}}(\tau)= -AZ_{m,A}(\tau)+V_{m+1} T_{m+1,m}^A L_{m,A}(\tau),
\end{equation}
where $L_{m,A}(\tau)={\widetilde E}_me^{(t-\tau) {\cal T}_{m,A}}E_m$.\\
Therefore, the error $e_{m,A}(\tau)=Z_A(\tau)-Z_{m,A}(\tau)$ is such that
$$e_{m,A}'(\tau)=-Ae_m(\tau)-V_{m+1} T_{m+1,m}^A L_{m,A}(\tau),$$
which gives the following expression of $e_m$:
\begin{equation}
e_{m,A}(\tau)= -\displaystyle \int_0^{\tau} e^{(u-\tau)A}V_{m+1} T_{m+1,m}^A L_{m,A}(u) d u.
\end{equation}
On the other hand, since $\tau-u >0$, it follows that 
$$\Vert e^{(u-\tau)A} \Vert \le  e^{(\tau-u)\mu_2(-A)}= e ^{(u-\tau)\nu(A)}.$$
Then, we get 
$$\Vert e_{m,A}(\tau) \Vert \le \Vert T_{m+1,m}^A \Vert \int_0^{\tau}  e ^{(u-\tau)\nu(A)} \Vert L_{m,A}(u) \Vert d u.$$

\noindent Notice that if  $\nu(A)$ is not known but $\nu(A) \ \ge 0$ (which is the case for positive semidefinite matrices)  then we get the upper bound
\begin{equation}\label{erm2}
\Vert e_{m,A}(\tau) \Vert \le \Vert T_{m+1,m}^A \Vert \int_0^{\tau}   \Vert L_{m,A}(u) \Vert d u.
\end{equation}
To define a new upper bound for the norm of the global error ${\cal E}_m(t)$,  we can use the  upper bounds for the errors $e_{m,A}$ and $e_{m,B}$ in the expression \eqref{err3} stated in Propostion \ref{tr4} to get 
\begin{eqnarray*}
\Vert {\cal E}_m(t) \Vert & \le &  \Vert F\,  \Vert e^{t\mu_2(B)} \, \displaystyle \int_{t_0}^t e^{-\tau\mu_2(B)} \,  \Vert e_{m,A}(\tau)\Vert  d \tau \\ & + &  \Vert E_m\,  \Vert e^{t\mu_2(A)} \displaystyle \int_{t_0}^t e^{-\tau\mu_2(A)}  \Vert e_{m,B}(\tau)\Vert   d \tau,
\end{eqnarray*}
and then we obtain  
\begin{eqnarray}
\label{globerr}
\Vert {\cal E}_m(t) \Vert  & \le &  \Vert F\,  \Vert e^{t\mu_2(B)} \,\Vert T_{m+1,m}^A \Vert  \displaystyle \int_{t_0}^t e^{-\tau\mu_2(B)} \, S_{m,A}(\tau)  d \tau \\ & + &  \Vert E_m\,  \Vert e^{t\mu_2(A)} \Vert T_{m+1,m}^B \Vert\displaystyle \int_{t_0}^t e^{-\tau\mu_2(A)} S_{m,B} (\tau)  d \tau,
\end{eqnarray}
where $S_{m,A}(\tau)=\displaystyle  \int_0^{\tau}  e ^{(u-\tau)\nu(A)} \Vert L_{m,A}(u) \Vert d u$ and  $S_{m,B}(\tau)= \displaystyle  \int_0^{\tau}  e ^{(u-\tau)\nu(B)} \Vert L_{m,B}(u) \Vert d u$.\\
As $m$ is generally very small as compared to $n$ and $p$, the factors $L_{m,A}$ and $L_{m,B}$ can be computed using Matlab fuctions such as {\tt expm} and the integral appearing in the right sides of \eqref{erm}  and  \eqref{globerr}, can be approximated via  a quadrature formulae.\\

We summarize the steps of our proposed first approach (using  the extended block Arnoldi) in the following algorithm
\begin{algorithm}[h!]
\caption{The  extended block Arnoldi (EBA-exp) method for DSE's}\label{algo_EBA_exp}
\begin{itemize}
\item Input $X_0=X(t_0)$, a tolerance $tol>0$, an integer $m_{max}$.
\item  For $ m = 1,\ldots,m_{max} $
\begin{itemize}
\item  Apply the extended block Arnoldi algorithm to $(A,E)$ and $(B^T,F)$ to get the orthonormal matrices ${\mathcal V}_m=[V_1,...,V_m]$ and ${\mathcal W}_m=[W_1,...,W_m]$ and the upper block Hessenberg matrices  ${\mathcal T}_{m,A}$ and ${\mathcal T}_{m,B}$.
\item Set ${E}_m={\mathcal V}_m^TE$, ${F}_m={\mathcal W}_m^TF$ and compute $ Z_{m,A}(\tau)= e^{(t-\tau)\mathcal{T}_{m,A}}E_m$ and $ Z_{m,B}(\tau)= e^{(t-\tau)\mathcal{T}_{m,B}}F_m$ using the matlab function {\tt expm}.
\item Use a quadrature method to compute the integral \eqref{gm} and get an approximation of $G_m(t)$ for each $t \in [t_0,\, T_f]$.
\item If $\parallel R_m(t) \parallel =  \displaystyle \max \{ \parallel T_{m+1,m}^A {\bar G}_{m}(t) \parallel,  \parallel T_{m+1,m}^B {\bar G}_{m}(t) \parallel\} < tol$ stop and compute the approximate solution $X_m(t)$ in the factored form given by the relation \eqref{approx}.
\end{itemize}
\item End
\end{itemize}
\end{algorithm}

\section{Projecting  and solving the low dimensional problem}

\subsection{Low-rank approximate solutions}\label{ss3.1}
In this section, we show how to obtain low rank approximate solutions to the differential Sylvester equation \eqref{sylv1} by first  projecting directly the initial problem onto   block (or extended block) Krylov subspaces and then solve the obtained low dimensional differential problem.  
We first apply the  block Arnoldi  algorithm (or the extended block Arnoldi)  to the pairs $(A,E)$ and $(B^T,F)$ to get the orthonormal matrices ${\cal V}_m$ and  ${\cal W}_m$,  whose columns form orthonormal bases of the extended block Krylov subspaces ${\cal K}_m(A,E)$ and  ${\cal K}_m(B^T,F)$, respectively. We also get the upper block Hessenberg matrices $ {\cal T}_{m,A}={\cal V}^T_m A {\cal V}_m $ and $ {\cal T}_{m,B}={\cal W}^T_mB^T {\cal W}_m $. \\ Let $X_m(t)$ be the desired low rank approximate solution  given as 
\begin{equation}\label{approx1}
X_m(t) = {\cal V}_m Y_m(t) {\cal W}_m^T,
\end{equation}
satisfying the Petrov-Galerkin orthogonality condition
\begin{equation}
\label{galerkin}
{\cal V}_m^T R_m(t) {\cal W}_m =0,\; t \in [t_0,\; T_f],
\end{equation}
where $R_m(t)$ is the residual $ R_m(t) = \displaystyle {\dot X}_m(t)-A\,X_m(t)-X_m(t)\,B- EF^T $.  Then, from \eqref{approx1} and \eqref{galerkin}, we obtain the low dimensional differential Sylvester  equation
\begin{equation}\label{lowsylv}
\displaystyle {\dot Y}_m(t)- {\cal T}_{m,A}\,Y_m(t)-Y_m(t)\,{\cal T}_{m,B}^T  - E_mF_m^T=0,
\end{equation}
 where  $ { E}_m= {\cal V}_m^T\,E$ and $ { F}_m= {\cal W}_m^T\,F$. The obtained low dimensional differential Sylvester equation \eqref{lowsylv} is  the  same as the one given by \eqref{low2}. We have now to solve the latter differential  equation by some integration method such as the well known Backward Differentiation Method (BDF) \cite{butcher}  or the  Rosenbrock method \cite{butcher,rosenbrock}. \\  Notice that all the properties and results such as the expressions of the residual norms or the upper bounds for the norm of the error given in the  last section are still valid with this second approach. The two approaches only differ in the way  the projected low dimensional differential Sylvester matrix equations are numerically solved.

\subsection{BDF for solving the low order differential Sylvester equation \eqref{lowsylv}}\label{projbdf}

We use the  Backward Differentiation Formula (BDF) method for solving, at each step $m$ of the  extended block  Arnoldi (or block Arnoldi) process,  the low dimensional differential Sylvester matrix equation \eqref{lowsylv}.  We notice that BDF is  especially used for the solution of stiff differential equations.\\  At each time $t_k$, let  $Y_{m,k}$ of the approximation of $Y_m(t_k)$, where $Y_m$ is a  solution of  (\ref{lowsylv}).  Then, the new approximation $Y_{m,k+1}$ of  $Y_m(t_{k+1})$ obtained at step $k+1$ by BDF is defined  by the implicit relation 
\begin{equation}
\label{bdf}
Y_{m,k+1} = \displaystyle \sum_{i=0}^{p-1} \alpha_i Y_{m,k-i} +h_k \beta {\mathcal F}(Y_{m,k+1}),
\end{equation} 
where $h_k=t_{k+1}-t_k$ is the step size, $\alpha_i$ and $\beta_i$ are the coefficients of the BDF method as listed  in Table \ref{tab1} and ${\mathcal F}(Y)$ is  given by 
$${\mathcal F}(Y)= {\cal T}_{m,A}\,Y+Y\,{\cal T}_{m,B}^T+\,E_m\,F_m^T.$$

\begin{table}[h!!]
\begin{center}
\begin{tabular}{c|cccc} 
\hline
$p$ & $\beta$ &$\alpha_0$ & $\alpha_1$ & $\alpha_2$ \\
\hline
1 & 1 & 1 & &\\
2 & 2/3 & 4/3& -1/3 &\\
3 & 6/11 & 18/11 & -9/11 & 2/11\\
\hline
\end{tabular}
\caption{Coefficients of the $p$-step BDF method with $p \le 3$.}\label{tab1}
\end{center}
\end{table}
\noindent The approximate $Y_{m,k+1}$ solves the following matrix equation
\begin{equation*}
-Y_{m,k+1} +h_k\beta ({\cal T}_{m,A} Y_{m,k+1} + Y_{m,k+1} {\cal T}_{m,B}^T+ EF^T) + \displaystyle \sum_{i=0}^{p-1} \alpha_i Y_{m,k-i} = 0,
\end{equation*}
which can be written as the following  Sylvester  matrix equation

\begin{equation}
\label{sylvbdf}
\mathbb{T}_{m,A}\, Y_{m,k+1}  + \,Y_{m,k+1} \mathbb{T}_{m,B}^T+ \mathbb{E}_{m,k}\, \mathbb{F}_{m,k}^T   =0.
\end{equation}
We assume  that at each time $t_k$, the approximation $Y_{m,k}$ is  factorized  as a low rank product 
 $Y_{m,k}\approx {\widetilde U}_{m,k} {{\widetilde V}_{m,k}}^T$, where ${\widetilde U }_{m,k} \in \mathbb{R}^{n \times m_k}$ and ${\widetilde V }_{m,k} \in \mathbb{R}^{p \times m_k}$, with $m_k \ll n,p$. In that case, the coefficient matrices appearing in \eqref{sylvbdf} are given by
$$\mathbb{T}_{m,A}= h_k\beta {\cal T}_{m,A} -\displaystyle \frac{1}{2}I;\;\;   \mathbb{T}_{m,B}= h_k\beta {\cal T}_{m,B} -\displaystyle \frac{1}{2}I,$$
$$ \mathbb{E}_{m,k+1}=[\sqrt{h_k\beta} E^T, \sqrt{\alpha_0}{\widetilde U }_{m,k}^T,\ldots,\sqrt{\alpha_{p-1}} {\widetilde U }_{m,k+1-p}^T]^T\, $$
 and 
$$\mathbb{F}_{m,k+1}=[\sqrt{h_k\beta} F^T, \sqrt{\alpha_0}{\widetilde V }_{m,k}^T,\ldots,\sqrt{\alpha_{p-1}} {\widetilde V }_{m,k+1-p}^T]^T.$$
The   Sylvester matrix  equation \eqref{sylvbdf} can be solved by applying direct methods based on Schur decomposition such as the Bartels-Stewart algorithm \cite{bartels,gnv}. \\ 
Notice that we can also use the BDF method applied directly to the original problem \eqref{sylv1} and then at  each iteration, one has to solve large Sylvester matrix equations which  can be done by using Krylov-based methods as developed in \cite{Elguen,jbilou}.

\subsection{Solving the low dimensional problem with the Rosenbrock method}\label{projros1}
Applying Rosenbrock method \cite{butcher,rosenbrock}  to the low dimensional differential Sylvester matrix equation \eqref{lowsylv},  the new approximation $Y_{m,k+1}$ of  $Y_m(t_{k+1})$ obtained at step $k+1$  is defined, in the ROS(2) particular case by the relations\\
\begin{equation}\label{ros1}
Y_{m,k+1} =Y_{m,k}+ \displaystyle \frac{3}{2}K_1+ \frac{1}{2}K_2,
\end{equation}
where $K_1$ and $K_2$ solve the following Sylvester equations
\begin{equation}\label{ros2}
\widetilde {\mathbb{T}}_{m,A}K_1+K_1\widetilde {\mathbb{T}}_{m,B}= -\mathcal{F}(t_k,Y_{m,k}),
\end{equation}
and
\begin{equation}\label{ros3}
\widetilde {\mathbb{T}}_{m,A}K_2+K_2\widetilde {\mathbb{T}}_{m,B}= -\mathcal{F}(t_{k+1},Y_{m,k}+K_1)+ \displaystyle \frac{2}{h}K_1,
\end{equation}
where
$$\widetilde {\mathbb{T}}_{m,A}= \gamma {\mathcal{T}}_{m,B}-\displaystyle \frac{1}{2h}I\;\;{\rm and} \;\;\widetilde {\mathbb{T}}_{m,B} = \gamma  {\mathcal{T}}_{m,B}^T-\displaystyle \frac{1}{2h}I,$$
and
$$\mathcal{F}(Y)= {\mathcal{T}}_{m,A}Y+Y {\mathcal{T}}_{m,B}^T+E_mF_m^T.$$

\noindent We summarize the steps of the second approach (using  the extended block Arnoldi) in the following algorithm
\begin{algorithm}[h!]
\caption{The  extended block Arnoldi (EBA) method  for DSE's}\label{algo_EBA}
\begin{itemize}
\item Input $X_0=X(t_0)$, a tolerance $tol>0$, an integer $m_{max}$.
\item  For $ m = 1,\ldots,m_{max} $
\begin{itemize}
\item  Apply the extended block Arnoldi algorithm to the pairs $(A,E)$ and $(B^T,F)$ to compute  the  orthonormal bases ${\mathcal V}_m=[V_1,...,V_m]$ and   ${\mathcal W}_m=[W_1,...,W_m]$ and also the  the upper block Hessenberg matrices ${\mathcal T}_{m,A}$ and  ${\mathcal T}_{m,B}$.
\item Use the BDF or the Rosenbrock method to solve the low dimensional differential Sylvester  equation
$$\displaystyle {\dot Y}_m(t)- {\cal T}_{m,A}\,Y_m(t)-Y_m(t)\,{\cal T}_{m,B}^T  - E_mF_m^T=0,\; t \in [t_0,\,T_f]$$
\item If $\parallel R_m(t) \parallel  < tol$ stop and compute the approximate solution $X_m(t)$ in the factored form given by the relation \eqref{approx}.
\end{itemize}
\item End
\end{itemize}
\end{algorithm}

\newpage
\section{Numerical examples}

In this section, we compare the  approaches presented in this paper. The exponential approach (EBA-exp) summarized in Algorithm \ref{algo_EBA_exp}, which is based on the approximation of the solution to \eqref{sylv1} applying a quadrature method to compute the projected exponential form solution \eqref{gm}.  We used a scaling and squaring strategy, implemented in the MATLAB  \textbf{expm} function; see  \cite{higham05,moler03} for more details.  The second method (Algorithm 2) is based on the BDF  integration method applied to the projected Sylvester equation as described in Section (\ref{projbdf}). Finally, we considered the EBA-ROS(2) method as described in Section (\ref{projros1}). The basis of the projection subspaces were generated by the extended block Arnoldi algorithm for all methods.
All the experiments were performed on a laptop with an  Intel Core i7 processor and 8GB of RAM. The algorithms were coded in Matlab R2014b. \\

\noindent {\bf Example 1}. 
For this example, the matrices  $A \in \mathbb{R}^{n \times n}$ and $B\in \mathbb{R}^{p \times p}$ were   obtained from the 5-point discretization of the operators 
\begin{equation*}
\displaystyle{L_A=\Delta u-f_1(x,y)\frac{\partial u}{\partial x}+ f_2(x,y)\frac{\partial u}{\partial y}+g_1(x,y)},
\end{equation*}
and
\begin{equation*}
\displaystyle{L_B=\Delta u-f_3(x,y)\frac{\partial u}{\partial x}+ f_4(x,y)\frac{\partial u}{\partial y}+g_2(x,y)},
\end{equation*}
 on the unit square $[0,1]\times [0,1]$ with homogeneous Dirichlet boundary conditions.  The number of inner grid points in each direction are  $n_0$ for $A$ and $p_0$ for $B$ and the dimension of the matrices $A$ and $B$ are $n = n_0^2=$ and $p=p_0^2$ respectively. Here we set $f_1(x,y) = x+10y^2$, $\displaystyle {f_2(x,y)= \sqrt{2x^2+y^2}}$, $f_3(x,y) = x+2y$, $f_4(x,y)= exp(y-x)$ ,  $g_1(x,y) = x^2-y^2$ and $g_2(x,y)=y^2-x^2$.  The time interval considered was $[0,\,2]$ and the initial condition $X_0=X(0)$ was  $X_0=Z_0Z_0^T$, where $Z_0=0_{n \times 2}$.\\

\begin{figure}[h!]
	\begin{center}
		\includegraphics[width=16cm,height=11cm]{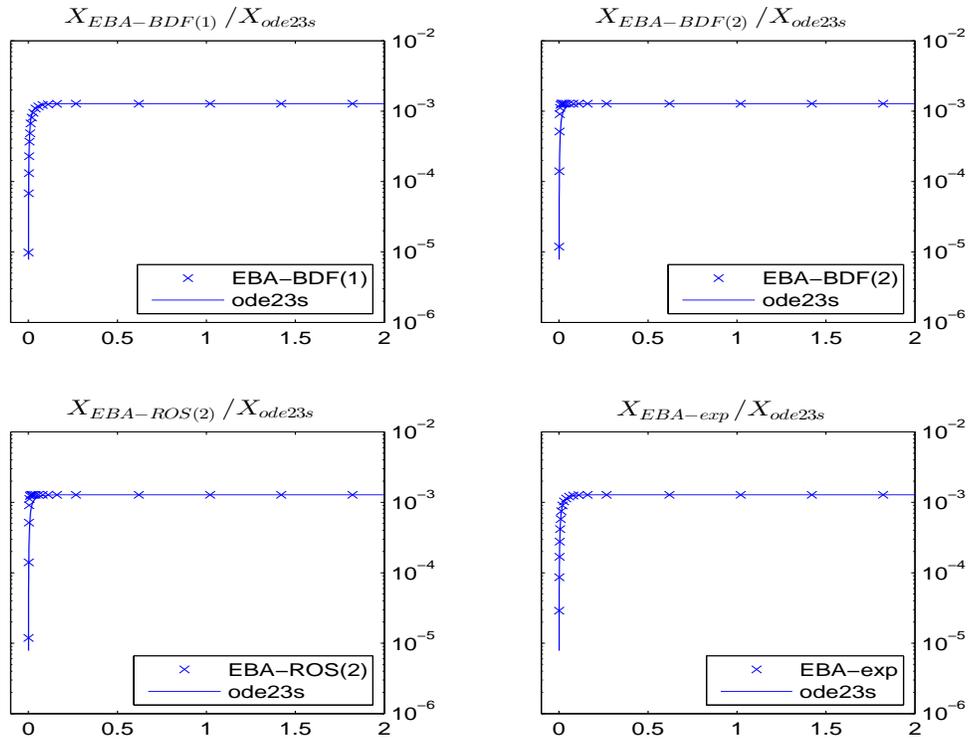}
		\caption{ Values of $X_{11}(t)$ for  $t \in [0,\, 2]$}\label{Figure2}
	\end{center}
\end{figure}

\noindent For all projection-based methods, we used  projections onto the Extended Block Krylov subspaces ${\mathcal K}_k(A,B) = {\rm Range}(B,A\,B,\ldots,A^{m-1}\,B,A^{-1}\,B,\ldots,(A^{-1})^m\,B)$ and the  tolerance was set to $10^{-10}$ for the stop test on the residual.  For the EBA-BDF and Rosenbrock methods, we used a constant timestep $h$. The entries of the matrices $E$ and $F$ were random values uniformly distributed on the interval $[0, \, 1]$ and their rank were set to $s=2$. \\
To the authors' knowledge, there are no available exact solutions of large scale matrix Sylvester differential equations in the
literature. In order to check if our approaches produce reliable results, we first compared our results to the one given by Matlab's ode23s solver which is designed for stiff differential equations. This was done by vectorizing our DSE, stacking the columns of $X$ one on top of each other. This method is not suited to large-scale problems. Due to the memory limitation of our computer when running the ode23s routine,  we chose a  size of $100\times 100$ for the matrices $A$ and $B$.\\
\noindent In Figure \ref{Figure2}, we compared the component $X_{11}$ of the solution obtained by the methods tested in this section, to the solution provided by the ode23s method from Matlab,  on the time interval $[0,\,2]$, for $size(A),~size(B)=100\times 100$ and a constant timestep $h=10^{-2}$. We observe that all the considered methods give good results in terms of accuracy. The relative error norms  at final time $T_f=2$ were of order $\mathcal{O} (10^{-10})$ for the EBA-exp method and $\mathcal{O} (10^{-12})$ for the others.  The runtimes were respectively 0.6s for EBA-exp, 7.3s for EBA-BDF(1), 20.8s for EBA-BDF(2) and 29.2s for EBA-ROS(2). The ode23s routine required 978s.
\noindent In Table \ref{tab2}, we give  the obtained runtimes in seconds, the number of Arnoldi iterations and the Frobenius residual norm at final time, for the resolution of Equation \eqref{sylv1}  for $t \in [0,\, 2]$, with a timestep $h=0.01$.  
\begin{table}[h!]
	\begin{center}
		\begin{tabular}{c | c c c c }
		 &EBA-exp&EBA-BDF(1)&EBA-BDF(2)&EBA-ROS(2)\\
			\hline
			$n,p=2500,2500$&$3.8$s $(m=16)$&$6.1$s $(m=18)$&$13.6$s $(m=18)$&$28.8$s $(m=23)$\\
			$\|R_m(T_f)\|_F$&$1.04\times 10^{-8}$ &$2.45\times 10^{-10}$&$2.45\times 10^{-10}$&$3.05\times 10^{-10}$\\
			\hline
				$n,p=10000,10000$&$35.2$s $(m=22)$&$38.4$s $(m=25)$&$80.3$s $(m=25)$&$104.7$s $(m=33)$\\
			$\|R_m(T_f)\|_F$&$4.4\times 10^{-9}$ &$4.1\times 10^{-11}$&$4.2\times 10^{-11}$&$5.8\times 10^{-11}$\\
			\hline
				$n,p=22500,10000$&$137.3$s $(m=22)$&$166.5$s $(m=30)$&$342.3$s $(m=30)$&$246$s $(m=35)$\\
			$\|R_m(T_f)\|_F$&$4.2\times 10^{-8}$ &$3.7\times 10^{-11}$&$3.6\times 10^{-11}$&$1.78\times 10^{-9}$\\
			\hline
		\end{tabular}
		\caption{Runtimes in seconds and the residual norms}\label{tab2}
	\end{center}
\end{table}
\noindent The results  in Table \ref{tab2} show that the EBA-exp method is outperformed by the other approaches in terms of accuracy, although it allows to obtain an acceptable approximation more quickly. The EBA-BDF(1) appears to be the better option in terms of time and accuracy.\\

\noindent {\bf Example 2} 
In this second example, we considered the particular case \begin{equation}\label{sylvrail}
\left\{
\begin{array}{l}
\dot X(t)=A(t)\,X(t)+X(t)\,A(t)-E(t)F(t)^T;\; (DSE) \\
\;X(t_0)=X_0,\; \; t \in [t_0, \, T_f],
\end{array}
\right.
\end{equation}
where the matrix  $A=Rail1357$ was extracted from  the IMTEK collection Optimal Cooling of Steel Profiles \footnote{https://portal.uni-freiburg.de/imteksimulation/downloads/benchmark}.  We compared the EBA-BDF(1) method to the EBA-exp and EBA-ROS(2) methods for the problem size $n=1357$  on the time interval $[0\,,2]$. The initial value $X_0$ was chosen as $X_0=0$ and  the timestep was set to $h=0.001$. The tolerance for EBA  stop test was set to $10^{-7}$ for all methods and the projected low dimensional Sylvester equations were numerically solved by the  solver ({\tt lyap} from Matlab at each iteration of the extended block Arnoldi  algorithm for the EBA-BDF(1), EBA-BDF(2) and EBA-ROS(2) methods. As the size of the coefficient matrices allowed it, we also computed an approximate solution of \eqref{sylvrail} applying a quadrature method to the integral form of the exact solution given by Formula\eqref{solexacte2} and took it as a reference solution. 
In Table \ref{tab3}, we reported the runtimes, in seconds, the number $m$ of Arnoldi iterations and the Frobenius norm $\|\mathcal{E}(T_f)_m\|_F$ of the error at final time.
\begin{table}[h!]
	\begin{center}
		\begin{tabular}{c | c c c c }
		&EBA-exp&EBA-BDF(1)&EBA-BDF(2)&EBA-ROS(2)\\
			\hline
			Runtime (s)&$48.4$ s ($m=18$)&$471.9$ s  ($m=18)$&$1549.2$s ($m=23$) & $1827$s ($m=21$)\\
			$\|\mathcal{E}_m(T_f)\|_F$&$1.28 \times 10^{-10}$& $5 \times 10^{-5}$ & $1.48 \times 10^{-4}$&$4.9 \times 10^{-5}$.\\
			\hline
		\end{tabular}
		\caption{Optimal Cooling of Steel Profiles: runtimes, number of Arnoldi iterations and error norms }\label{tab3}
	\end{center}
\end{table}
\noindent As can be seen from the reported  results in Table \ref{tab3}, the EBA-exp method clearly outperforms all the other listed options. 

\begin{figure}[h!]
\begin{center}
\includegraphics[width=6cm,height=4cm]{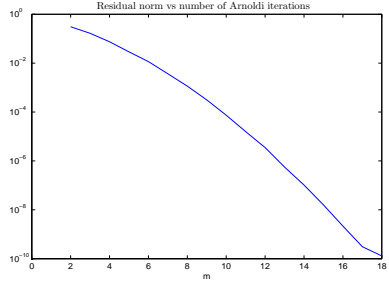}
\caption{ Residual norm \textit{vs} number $m$ of Arnoldi iterations}\label{Figure3}
\end{center}
\end{figure}
\noindent In Figure \eqref{Figure3}, we plotted the Frobenius residual norm $\|R_m(t_f)\|_F$ at final time $T_f$ in function of the number $m$ of Arnoldi iterations for the EBA-exp method.
%\newpage

\section{Appendix A}

Here we recall the extended block Arnoldi (EBA) and block Arnoldi (BA) algorithms,  when applied to the pair $(A,E)$. EBA is described in Algorithm \ref{eba} as follows
\begin{algorithm}[h!]
\caption{The extended block Arnoldi algorithm (EBA)}\label{eba}
\begin{itemize}
\item   Inputs: $ A $ an $ n \times n $ matrix, $ E $ an $ n \times s $ matrix and $ m $ an integer.
\item Compute the QR decomposition of $ [E,A^{-1}E] $,\textit{ i.e}., $ [E,A^{-1}E] = V_1\Lambda $; \\
\hspace*{1.2cm} Set $ {\mathcal V}_0 = \left [ ~\right] $;
\item For $ j = 1,\ldots,m $
%\hspace*{1.8cm} $ {\ cal V}_j = \[{\mathcal V}_{j-1}, V_j\] $; \\
\item \hspace*{0.4cm} Set $ V_j^{(1)} $: first $ s $ columns of $ V_j $ and  $ V_j^{(2)} $: second $ s $ columns of $ V_j $
\item  \hspace*{0.4cm} $ {\mathcal V}_j = \left [ {\mathcal V}_{j-1}, V_j \right ] $; $ \hat V_{j+1} = \left [ A\,V_j^{(1)},A^{-1}\,V_j^{(2)} \right ] $.
\item  Orthogonalize $ \hat V_{j+1} $ w.r.t $ {\mathcal V}_j $ to get $ V_{j+1} $, \textit{i.e.},\\
\hspace*{1.8cm} For $ i=1,2,\ldots,j $ \\
\hspace*{2.4cm} $ H_{i,j} = V_i^T\,\hat V_{j+1} $; \\
\hspace*{2.4cm} $ \hat V_{j+1} = \hat V_{j+1} - V_i\,H_{i,j} $; \\
\hspace*{1.8cm} Endfor $i$
\item Compute the QR decomposition of $ \hat V_{j+1} $, \textit{i.e.}, $ \hat V_{j+1} = V_{j+1}\,H_{j+1,j} $.\\
\item Endfor $j$.
\end{itemize} ${}$ 
\end{algorithm}
\vskip0.3cm

The block Arnoldi algorithm is summarized in Algorithm \ref{ba} as follows
 \begin{algorithm}[h!]
 \caption{The block Arnoldi algorithm (BA)}\label{ba}
 \begin{itemize}
 \item   Inputs: $ A $ an $ n \times n $ matrix, $ E $ an $ n \times s $ matrix and $ m $ an integer. 
 \item Compute the QR decomposition of $ E $,\textit{ i.e}., $ E= V_1 R_1 $.
\item For $j=1,\ldots,m$
\begin{enumerate}
\item $W_j=AV_j$,
 \item for $i=1,2,\ldots,j$
 \begin{itemize}
\item $H_{i,j}=V_i^T\,W_j,$
\item $W_j=W_j-V_j\,H_{i,j},$
\end{itemize}
\item endfor 
\item $Q_j R_j=W_j$ ($QR$ decomposition)
\item $V_{j+1}=Q_j$, and  $H_{j+1,j}=R_j$.
\end{enumerate}
\item EndFor $j$
\end{itemize}
\end{algorithm}
\medskip

\noindent Since the above algorithms implicitly involve a Gram-Schmidt process, the obtained blocks $ {\cal V}_m = \left [
V_1,V_2,\ldots,V_m \right ] $ ($ V_ i \in \mathbb{R}^{n \times d} $) ,where $d=s$ for the block Arnoldi and $d=2s$ for the extended block Arnoldi, have their columns mutually orthogonal provided none of the upper triangular  matrices $ H_{j+1,j} $ are rank deficient. 
Hence, after $ m $ steps, Algorithm \ref{eba} and Algorithm \ref{ba}  build orthonormal bases $ {\cal V}_m $ of the Krylov subspaces ${\cal K}_m(A,E)={\rm Range}(E,A\,E,\ldots,A^{m-1}\,E,A^{-1}\,E,\ldots,(A^{-1})^m\,E) $ or  $ \K_{m}(A,E) ={\rm Range}(E,A\,E,\ldots,A^{m-1}\,E) $, respectively  and a block upper Hessenberg matrix $ {\cal H}_m $ whose nonzero sub-blocks are the $ H_{i,j} $. Note that each submatrix $ H_{i,j} $ ($1 \le i \le j \le m $) is of order $ d$. \\
\noindent Let $ {\cal T}_m \in \mathbb{R}^{d \times d} $ be the
restriction of the matrix $ A $ to the extended Krylov  subspace $
{\cal K}_m(A,E) $ (or to the block Krylov subspace $\K_{m}(A,E)$), i.e., $ {\cal T}_m = {\cal V}_m^T\,A\,{\cal
V}_m $. Then it can be shown that  matrix  $ {\cal T}_m $ is
also block upper Hessenberg with $ d \times d $ blocks, see\cite{heyouni09,simoncini1}  .
For the block Arnoldi algorithm, ${\cal T}_m={\cal H}_m$  while for the extended block Arnoldi algorithm,  a recursion can be  derived to compute $ {\cal T}_m $ from $
{\cal H}_m $ without requiring matrix-vector products with $A $, see \cite{simoncini1}. We notice that for large and non structured problems, the inverse of the matrix $ A $ is not computed explicitly and in this
case we can use iterative solvers with preconditioners to solve
linear systems with $ A $.

\section{Conclusion}
We presented in the present paper two new approaches for computing approximate solutions to large scale differential  Sylvester  matrix equations. The first one comes naturally from the exponential expression of the exact solution and the use of approximation techniques of the exponential of a matrix times a block of vectors. The second approach is obtained by first projecting the initial problem onto a block Krylov (or extended Krylov) subspace, obtain a low dimensional differential Sylvester  equation which is solved by using the well known BDF or Rosenbrock integration method. We gave some theoretical results such as the exact expression of the residual norm and also upper bounds for the norm of the error.  Numerical experiments show that both approaches are promising for large-scale problems, with a clear advantage for the EBA-exp method in terms of computation time although the EBA-BDF(1) method shows to offer a good balance between the execution time and the accuracy in some cases.

\bibliographystyle{plain}
\end{document}